\documentclass[12pt]{amsart}
\usepackage{amssymb,amsfonts,amsmath,amsopn,amstext,amscd,color,
latexsym,mathrsfs,skak,stackrel,stmaryrd,url,verbatim}
\usepackage{mathabx}
\usepackage[utf8]{inputenc}

\usepackage[all]{xy}

\renewcommand{\labelenumi}{{\it (\roman{enumi})}}
\renewcommand{\theenumi}{{\it (\roman{enumi})}}

\theoremstyle{remark}

\newtheorem{para}{\bf}[subsection]

\newtheorem{rems}[para]{\bf Remarks}
\newtheorem{rem}[para]{\bf Remark}
\newtheorem{convention}[para]{\bf Convention}

\theoremstyle{definition}

\newtheorem{dfn}[para]{Definition}

\theoremstyle{plain}

\newtheorem{thm}[para]{Theorem}

\newtheorem{lemma}[para]{Lemma}
\newtheorem{cor}[para]{Corollary}
\newtheorem{prop}[para]{Proposition}

\newenvironment{numequation}{\addtocounter{para}{1}
\begin{equation}}{\end{equation}}

\newcommand{\vpi}{{\varpi}}

\newcommand{\bbF}{{\mathbb F}}
\newcommand{\bbG}{{\mathbb G}}

\newcommand{\bbN}{{\mathbb N}}

\newcommand{\bbP}{{\mathbb P}}
\newcommand{\bbQ}{{\mathbb Q}}

\newcommand{\bbX}{{\mathbb X}}

\newcommand{\bbZ}{{\mathbb Z}}

\newcommand{\bare}{{\bar{e}}}
\newcommand{\baru}{{\bar{u}}}
\newcommand{\barx}{{\bar{x}}}

\renewcommand{\frm}{{\mathfrak m}}

\newcommand{\frp}{{\mathfrak p}}

\newcommand{\cC}{{\mathcal C}}

\newcommand{\cK}{{\mathcal K}}

\newcommand{\cO}{{\mathcal O}}
\newcommand{\cP}{{\mathcal P}}

\newcommand{\tA}{\tilde{A}}

\newcommand{\tK}{\tilde{K}}
\newcommand{\tk}{\tilde{k}}

\newcommand{\bcO}{\breve{\cO}}

\newcommand{\bQp}{{\breve{\bbQ}_p}}

\newcommand{\brF}{\breve F}

\newcommand{\Aut}{{\rm Aut}}

\newcommand{\car}{\stackrel{\simeq}{\longrightarrow}}

\newcommand{\Def}{{\rm Def}}

\newcommand{\End}{{\rm End}}

\newcommand{\Frac}{{\rm Frac}}

\newcommand{\Fp}{{\mathbb F_p}}
\newcommand{\Fpb}{{\overline{\mathbb F}_p}}

\newcommand{\Gal}{{\rm Gal}}
\newcommand{\GL}{{\rm GL}}

\newcommand{\Hom}{{\rm Hom}}
\newcommand{\hra}{\hookrightarrow}
\newcommand{\id}{{\rm id}}

\newcommand{\kb}{\overline{k}}

\newcommand{\lra}{\longrightarrow}

\newcommand{\midc}{{\; | \;}}
\newcommand{\ob}{{\rm ob}}

\newcommand{\Pf}{{\it Proof. }}

\newcommand{\prim}{{\rm prim}}

\newcommand{\Qp}{{\mathbb Q_p}}

\newcommand{\ra}{\rightarrow}

\newcommand{\sep}{{\rm sep}}

\newcommand{\Spf}{{\rm Spf}}

\newcommand{\univ}{{\rm univ}}

\newcommand{\Z}{{\mathbb Z}}

\newcommand{\Zp}{{\mathbb Z_p}}

\newcommand{\fm}{\mathfrak{m}}

\newcommand{\df}[2][]{\emph{#2}}

\newcommand{\llrrparen}[1]{
  \left(\mkern-4mu\left(#1\right)\mkern-4mu\right)}

\begin{document}

\title{Lubin-Tate Deformation Spaces and Fields of Norms}
\author{Annie Carter}
\address{University of Hawai`i at M\=anoa, Department of Mathematics, 2565 McCarthy Mall (Keller Hall 401A), Honolulu, HI 96822, U.S.A.}
\email{atcarter@hawaii.edu}
\author{Matthias Strauch}
\address{Indiana University, Department of Mathematics, Rawles Hall, Bloomington, IN 47405, U.S.A.}
\email{mstrauch@indiana.edu}

\thanks{A.C. would like to acknowledge support by NSF RTG grant DMS-1502651.}

\begin{abstract} We construct a tower of fields from the rings $R_n$ which parametrize pairs $(X,\lambda)$, where $X$ is a deformation of a fixed one-dimensional formal group $\bbX$ of finite height $h$, together with a Drinfeld level-$n$ structure $\lambda$. We choose principal prime ideals $\frp_n \midc (p)$ in each ring $R_n$ in a compatible way and consider the  field $K'_n$ obtained by localizing $R_n$ at $\frp_n$, completing, and passing to the fraction field. By taking the compositum $K_n = K'_n K_0$ of each field with the completion $K_0$ of a certain unramified extension of $K'_0$, we obtain a tower of fields $(K_n)_n$ which we prove to be strictly deeply ramified in the sense of Anthony Scholl. When $h=2$ we also investigate the question of whether this is a Kummer tower. 
\end{abstract}

\maketitle

\tableofcontents

\section{Introduction}
\setcounter{subsection}{1}

In this paper we study a tower $K_\bullet$ of complete discrete valuation fields of characteristic zero and residue field of characteristic $p>0$:

\begin{numequation}\label{tower_intro} K_0 \subseteq K_1 \subseteq K_2 \subseteq K_3 \subseteq \cdots 
\end{numequation}

The fields $K_n$ are defined in terms of torsion points of the universal deformation of a formal group $\bbX$ of dimension one and height $h>0$ over $\Fpb$. The extension $K_n/K_0$ is Galois with $\Gal(K_n|K_0) \simeq (\Z/p^n)^\times \ltimes (\Z/p^n)^{h-1}$. When $h=1$ we have $\bbX = \widehat{\bbG}_m$, and $K_n = \bQp(\mu_{p^n})$, where $\bQp = W(\Fpb)\left[\frac{1}{p}\right]$ is the completion of the maximal unramified extension of $\Qp$. If $h>1$ the residue field of $K_0$ is a separable (infinite) Galois extension of $\Fpb\llrrparen{u_1, \ldots, u_{h-1}}$ and thus imperfect. 

\vskip8pt

Our starting point is the paper \cite{Scholl} of A.J. Scholl in which he develops a theory of norm fields for certain towers of fields $(K_n)_n$ which he calls {\it strictly deeply ramified}\footnote{We refer to the body of the paper for the discussion of this concept.}, and whose residue fields are not necessarily perfect. Our first result is:

\vskip8pt

\begin{thm}\label{result1} The tower \ref{tower_intro} is strictly deeply ramified in the sense of A.J. Scholl. 
\end{thm}

By the general theory of \cite{Scholl}, the tower \ref{tower_intro} therefore gives rise to a complete discretely valued field of norms $E$ of characteristic $p$ whose ring of integers is

\[\cO_E = \varprojlim_n \cO_{K_n}/p\cO_{K_n} \;,\]

where the transition maps are the $p$-power maps.

\vskip8pt

It is an important question whether the field of norms $E$ lifts to characteristic zero in a way that is compatible with the action of $\Gamma = \Gal(K_\infty|K_0)$ and Frobenius (where, as usual, $K_\infty = \bigcup_n K_n$). If that is the case, then one can describe $p$-adic representations of $\Gal(\overline{K_0}|K_0)$ in terms of $(\phi,\Gamma)$-modules.

\vskip8pt

While we do not settle here the question whether our norm field $E$ lifts to characteristic zero (in a way that is compatible with the action of $\Gamma$ and Frobenius), we investigate if $K_\bullet$ is a Kummer tower. By this we mean that there are elements $t_1, \ldots, t_{h-1} \in K_0$ such that for all $n \ge 0$ one has $K_n = K_0(\mu_{p^n},\sqrt[p^n]{t_1}, \ldots, \sqrt[p^n]{t_{h-1}})$. Scholl shows in \cite[sec. 2.3]{Scholl}, that the field of norms associated to a Kummer tower lifts to characteristic zero.  

\vskip8pt

In order to explain the result that we have in this direction, we need to briefly 
sketch the construction of the tower \ref{tower_intro}; cf. section 
\ref{sec:construction-of-towers} for more details. Let $R_n$ be the ring 
which represents (isomorphism classes of) triples $(X,\iota,\lambda)$, 
where $(X,\iota)$ is a deformation of $\bbX$ and 
$\lambda: p^{-n}\Z/\Z \ra X[p^n]$ is a Drinfeld level-$n$ structure. The 
ring $R_0$ is a non-canonically isomorphic to 
$W(\Fpb)\llbracket u_1, \ldots, u_{h-1} \rrbracket$, and 
$R_n\left[\frac{1}{p}\right] / R_0\left[\frac{1}{p}\right]$ is an 
unramified Galois extension of rings with Galois group isomorphic to 
$\GL_n(\Z/p^n)$. Set $\frp_0 = p R_0$. For $n>0$ let $\frp_n \in R_n$ be 
a prime ideal of height one which has the property that $\frp_{n+1}$ 
divides $\frp_n$ in $R_{n+1}$, for all $n \ge 0$. Let $R'_n$ be the 
completion of the localization $(R_n)_{\frp_n}$ (with respect to the 
topology defined by the maximal ideal), and $K'_n = \Frac(R'_n)$ the 
field of fractions of $R'_n$. Let $K'_{n,u} \subseteq K'_n$ be the largest 
subfield which is unramified over $K'_0$, and put $\tK_0 = \bigcup_n K'_{n,u}$. Then define $K_0$ to be the $p$-adic completion of $\tK_0$ and define $K_n = K'_n K_0$ to be the composition of $K'_n$ and $K_0$. We have only investigated the question whether $K_\bullet = (K_n)_n$ is a Kummer tower when $h=2$. In this case, we have

\begin{thm}
\label{result2}
Let $h=2$ and let $K_\bullet = (K_n)_n$ be the tower as defined above. 
\vskip8pt

{\rm (i)} For every $n \geq 0$ there is an element $t_n \in K_0$ such that $K_n = K_0 \left(\mu_{p^n},\sqrt[p^n]{t_n}\right)$. \label{result2i}

\vskip5pt

{\rm (ii)} If $p>3$, there does not exist an element $t \in \tK_0$ such that for all sufficiently large $n$ one has $K_n = K_0\left(\mu_{p^n},\sqrt[p^n]{t}\right)$. \label{result2ii}
\end{thm}

\vskip8pt

Our method of proof leaves open the possibility that there is an element $t \in K_0$ such that $K_n = K_0\left(\mu_{p^n},\sqrt[p^n]{t}\right)$ for all $n \geq 0$. However, even if $K_\bullet$ fails to be a Kummer tower, it might still be possible that the norm field $E$ lifts to characteristic zero (together with Galois action and Frobenius), but we do not have positive evidence with regard to this problem.

\begin{rem} The tower of fields $(K'_n)_n$ and the field $K'_\infty = \bigcup_{n \ge 0} K'_n$ have also been studied in \cite{Kohlhaase} (where $K'_n$ is denoted $L_n$), but the objective in loc.\ cit.\ is quite different in that it is concerned with the composite field 

$$K'_\infty .K'_0\Big(\mu_{p^\infty},u_1^{1/p^\infty}, \ldots, u_{h-1}^{1/p^\infty}\Big) \;.$$

\vskip8pt

In particular, the fields $K_n$ do not appear in loc.\ cit., and the question whether the tower $K_\bullet$ is strictly deeply ramified, is not addressed in there.

\end{rem}

\section{Deformations with level structures and the tower of fields $K_\bullet$}\label{tower}

\subsection{Formal modules, deformations, and level structures} We briefly recall some facts about deformations of formal $\cO$-modules and level structures, following \cite[sections 1 and 4]{Drinfeld}. 

\vskip8pt
 
We fix a finite extension $F/\Qp$ with ring of integers $\cO = \cO_F$, uniformizer $\pi$, and residue field $k_F = \cO/(\pi)$ of cardinality $q$. Let $\bcO$ be the completion of the maximal unramified extension of $\cO$. The residue field of $\bcO$ is an algebraic closure of $k_F$ which we denote by $\kb_F$. We also fix a formal $\cO$-module $\bbX$ of dimension one and finite $F$-height $h \ge 1$ over $\kb_F$. Up to isomorphism there is only one formal $\cO$-module over $\kb_F$ of given $F$-height $h$ \cite[1.7]{Drinfeld}. Given a formal $\cO$-module $X$ over some ring $R$ we denote by $[\cdot]_X: \cO \ra \End_R(X)$ the corresponding ring homomorphism.

\vskip8pt

By $\cC$ we denote the category of $\bcO$-algebras $R$ with the following properties:

\begin{enumerate}\renewcommand{\theenumi}{\roman{enumi}}
	\renewcommand{\labelenumi}{(\theenumi)}
	\item $R$ is a complete, local, noetherian % needed?
	ring, whose maximal ideal we denote by $\fm_R$;
	\item the structure homomorphism $\bcO \to R$ is local; 
	\item the canonical field homomorphism $\kb_F = \bcO/\pi\bcO \to R/\fm_R$ is an isomorphism.
	\end{enumerate}

Morphisms in $\cC$ are local homomorphisms of $\bcO$-algebras.

\vskip8pt

By a \df{deformation} of $\bbX$ over $R \in \ob(\cC)$, we mean a pair $(X,\iota)$ consisting of a formal $\cO$-module $X$ over $R$, together with an isomorphism $\iota: \bbX \car X \otimes_R R/\frm_R$. Two deformations $(X_1,\iota_1)$ and $(X_2,\iota_2)$ are defined to be \df{isomorphic} if there is an isomorphism $f: X_1 \ra X_2$ of formal $\cO$-modules over $R$ such that $(f \otimes R/\frm_R) \circ \iota_1 = \iota_2$. In that case we write $f: (X_1,\iota_1) \car (X_2,\iota_2)$. 

\vskip8pt

Let $(X,\iota)$ be a deformation of $\bbX$. We fix a cooordinate $T$ on $X$, and using $T$, we equip the maximal ideal $\frm_R$ with the structure of an $\cO$-module. Let $n$ denote a positive integer. A \df[structure of level n@structure of level $n$]{structure of level $n$} on $X$ is an $\cO$-module homomorphism 
\[\lambda: \left(\pi^{-n}\cO/\cO\right)^h \lra \frm_R\] 

such that the power series $[\pi]_X(T)$ is divisible by

$$\prod_{\alpha \in \left(\pi^{-1}\cO/\cO\right)^h} (T-\lambda(\alpha)) \;.$$

\vskip8pt 

\begin{rem} A structure of level zero is, by definition, the unique homomorphism from the trivial group $(\pi^0\cO/\cO)^h$ to $\frm_R$. In the definition of $\Def_{\bbX,n}$ below, the datum of $\lambda$ can be ignored when $n=0$.
\end{rem}

Let $(X_1,\iota_1)$ and $(X_2,\iota_2)$ be two deformations of $\bbX$, and let $\lambda_i$ be a level-$n$ structure on $X_i$ for $i=1,2$. The triples $(X_1,\iota_1,\lambda_1)$ and $(X_2,\iota_2,\lambda_2)$ are defined to be \df{isomorphic} if there  is an isomorphism $f: (X_1,\iota_1) \car (X_2,\iota_2)$ of deformations satisfying $f \circ \lambda_1 = \lambda_2$. Define the functor 

\[\Def_{\bbX,n}: \cC \lra {\rm Sets}\]

\vskip8pt

by associating to $R \in \ob(\cC)$ the set of isomorphism classes of triples $(X,\iota,\lambda)$, where $(X,\iota)$ is a deformation of $\bbX$ and $\lambda$ is a level-$n$ structure on $X$. For $n' \ge n$, the restriction of any level-$n'$ structure $\lambda'$ on $X$ to $\left(\pi^{-n}\cO/\cO\right)^h \subseteq \left(\pi^{-n'}\cO/\cO\right)^h$ is a level-$n$ structure. We thus get a natural transformation $\Def_{\bbX,n'} \ra \Def_{\bbX,n}$. Moreover, we have a right action of $\GL_h(\cO/(\pi^n))$ on the functor $\Def_{\bbX,n}$ which is defined by $[X,\iota,\lambda].g = [X,\iota,\lambda \circ g]$, where $[X,\iota,\lambda]$ denotes the isomorphism class of the triple $(X,\iota,\lambda)$.

\vskip8pt

Parts (i)-(iii) of the following result are due to V.G. Drinfeld \cite[4.2, 4.3]{Drinfeld}, and part (iv) has been shown in \cite[2.1.2]{StrauchDeformationSpaces}.
 
\begin{thm}\label{thm:universal-deformation} {\rm (i)} For every $n \ge 0$ the functor $\Def_{\bbX,n}$ is representable, i.e., there is an $\bcO$-algebra $R_n \in \ob(\cC)$ and an isomorphism of functors 

$$\Def_{\bbX,n} \car \Hom_\cC(R_n, -) \;.$$

\vskip5pt

{\rm (ii)} The ring $R_n$ in (i) is a regular local ring. For all $n' \ge n$ the ring homomorphism  $R_n \ra R_{n'}$ (induced by the natural transformation $\Def_{\bbX,n'} \ra \Def_{\bbX,n}$) is finite and flat.

\vskip5pt

{\rm (iii)} The ring $R_0$ is (non-canonically) isomorphic to $\bcO\llbracket u_1, \ldots, u_{h-1}\rrbracket$.

\vskip5pt

{\rm (iv)} The ring extension $R_n\left[\frac{1}{\pi}\right] / R_0\left[\frac{1}{\pi}\right]$ is Galois with Galois group isomorphic to \linebreak $\GL_h(\cO/(\pi^n))$. (The left action of this group on $R_n$ is induced by its right action on the functor $\Def_{\bbX,n}$.) \qed
\end{thm}

\begin{rems} (i) When $F = \Qp$, hence $\cO = \Zp$, part (iii) is due to Lubin and Tate \cite{LubinTate}, which is why the formal scheme $\Spf(R_0)$ (or its rigid analytic generic fiber) is called a Lubin-Tate deformation space. More generally, the formal schemes $\Spf(R_n)$ (or their rigid analytic generic fibers) are also called Lubin-Tate deformation spaces. 

\vskip5pt

(ii) Let $[X^\univ,\iota^\univ] \in \Def_{\bbX,0}(R_0)$ be the element which corresponds to the identity map $\id_{R_0} \in \Hom_\cC(R_0,R_0)$. Then $X^\univ$ is called the \df{universal deformation} of $\bbX$. Furthermore, consider the isomorphism class of triples $[X^\univ,\iota^\univ,\lambda_n^\univ] \in \Def_{\bbX,n}(R_n)$ which corresponds to the identity map $\id_{R_n} \in \Hom_\cC(R_n,R_n)$. The map $\lambda_n^\univ$ is called the \df{universal level-$n$ structure}. Moreover, for $n' \ge n$, the restriction of $\lambda_{n'}^\univ$ to $\left(\pi^{-n}\cO/\cO\right)^h \subseteq \left(\pi^{-n'}\cO/\cO\right)^h$ is equal to the composition of $\lambda_n^\univ$ with the inclusion $\frm_{R_n} \hra \frm_{R_{n'}}$.

\vskip5pt

(iii) In the following we will often consider the action of $\GL_h(\cO)$ on $R_n$ which is induced by the canonical map $\GL_h(\cO) \ra \GL_h(\cO/(\pi^n))$, and we write $g.a$ for the image of $a \in R_n$ under the action of $g \in \GL_h(\cO)$, and we write $g.A$ for the image of a subset $A \subseteq R_n$ under the action of $g$.

\vskip5pt

(iv) When $h=1$ the universal deformation $X^\univ$ is, up to isomorphism, the unique lift of $\bbX$ to $\bcO$. This implies the well-known fact that all Lubin-Tate formal groups for $\cO$ (i.e., one-dimensional $\cO$-modules of $F$-height one over $\cO$) become isomorphic over $\bcO$, cf. \cite[sec. 3.7, Lemma 1]{Serre_LCFT}. Let $LT_\cO$ be any Lubin-Tate formal group for $\cO$, and let $F_n$ be the extension of $F$ generated by the $\pi^n$-torsion points of $LT_\cO$. This is a purely ramified extension of degree $(q-1)q^{n-1}$, and the composite field $\brF_n := F_n.\brF$, where $\brF = \bcO[\frac{1}{\pi}]$, does not depend on the choice of $LT_\cO$. When $h=1$, the ring $R_n$ is the ring of integers of $\brF_n$. \label{LT_extensions}
\end{rems}

\subsection{Construction of the tower of fields $K_\bullet$}\label{sec:construction-of-towers}

\begin{para} {\it Sequences of prime ideals.}\label{prime_ideals} The construction which we are going to perform depends on the choice of a sequence $\frp_\bullet  = (\frp_n)_{n > 0}$ of ideals $\frp_n \subseteq R_n$ with the following properties

\vskip8pt

\begin{enumerate}
\renewcommand{\theenumi}{\roman{enumi}}
\renewcommand{\labelenumi}{(\theenumi)}

\item For all $n>0$ the ideal $\frp_n$ is a prime ideal of height one.

\vskip5pt

\item $\frp_1 | (\pi)$ in  $R_1$, and $\frp_{n+1} | \frp_n$ in $R_{n+1}$ for all $n>0$.

\end{enumerate}
\end{para}

\vskip8pt

In the following we set $\frp_0 := \pi R_0$. We note that any prime ideal of height one of $R_n$ is a principal ideal, because $R_n$ is a regular local ring, hence a unique factorization domain. Put $R_\infty = \bigcup_{n \ge 0} R_n$.

\vskip8pt

\begin{para} Note that the group $\GL_h(\cO)$ acts on the set of all such sequences $\frp_\bullet$: if $g \in \GL_h(\cO)$, and if $\frp_\bullet = (\frp_n)_n$ is such a sequence, then $g.\frp_\bullet = (g.\frp_n)_n$ is another such sequence.
We call $\alpha  = (\alpha_1, \ldots, \alpha_h) \in (\pi^{-n}\cO/\cO)^h$ \df{primitive} if $\alpha$ is not divisible by $\pi$, i.e., $\alpha \not \in (\pi^{-(n-1)}\cO/\cO)^h$. Denote by $(\pi^{-n}\cO/\cO)^h_\prim \subseteq (\pi^{-n}\cO/\cO)^h$ the set of primitive elements. Note that the group of units $(\cO/(\pi^n))^\times$ acts on $(\pi^{-n}\cO/\cO)^h_\prim$, and let $\cP_n = (\pi^{-n}\cO/\cO)^h_\prim/(\cO/(\pi^n))^\times$ be the set of orbits under this group. We denote by $[\alpha] \in \cP_n$ the orbit of $\alpha \in (\pi^{-n}\cO/\cO)^h_\prim$.
We also call $v = (v_1, \ldots, v_h) \in \cO^h$ \df{primitive} if it is not divisible by $\pi$, we let $\cO^h_\prim$ be the subset of primitive vectors, and denote, as usual, by $\bbP^{h-1}(\cO) = \cO^h_\prim/\cO^\times$ the set of orbits under the action of $\cO^\times$, and we denote by $[v]$ its class in $\bbP^{h-1}(\cO)$.
\end{para}

Most statements of the following proposition have already been shown in the literature, but as we use them later on, we repeat  them here. For elements $x,y \in R_n$ we write $x \sim y$ if $x$ and $y$ are associate, i.e., there is $u \in R_n^\times$ such that $y = ux$ 

\begin{prop}\label{basic} Let $n$ be a positive integer. {\rm (i)} Let $\alpha_1, \ldots, \alpha_h \in (\pi^{-n}\cO/\cO)^h$ be a basis of $(\pi^{-n}\cO/\cO)^h$ over $\cO/(\pi^n)$. Then

$$\Big(\lambda_n^\univ(\alpha_1), \ldots, \lambda_n^\univ(\alpha_h)\Big)$$

\vskip8pt

is a regular system of parameters for $R_n$. Moreover, $R_n$ is generated as $R_0$-algebra by $\lambda_n^\univ(\alpha_1), \linebreak \ldots, \lambda_n^\univ(\alpha_h)$. 

\vskip5pt

{\rm (ii)} For every $\alpha \in (\pi^{-n}\cO/\cO)^h_\prim$, the element $\lambda_n^\univ(\alpha)$ is a prime element of $R_n$. 

\vskip5pt

{\rm (iii)} For $0 \neq \alpha \in (\pi^{-n}\cO/\cO)^h$ and $a \in (\cO/(\pi^n))^\times$, one has $\lambda_n^\univ(a\alpha) \sim \lambda_n^\univ(\alpha)$. Moreover, for  $\alpha, \beta \in (\pi^{-n}\cO/\cO)^h_\prim$ one has $\lambda_n^\univ(\alpha) \sim \lambda_n^\univ(\beta)$ if and only if $[\alpha] = [\beta]$ in $\cP_n$.

\vskip5pt

{\rm (iv)} For $\beta \in (\pi^{-n}\cO/\cO)^h_\prim$, and $\alpha \in (\pi^{-(n+1)}\cO/\cO)^h_\prim$, the prime element $\lambda^\univ_{n+1}(\alpha)$ divides $\lambda^\univ_n(\beta)$ in $R_{n+1}$ if and only if $[\pi\alpha] = [\beta]$ in $\cP_n$.

\vskip5pt

{\rm (v)} Let $\brF_n/\brF$ be as in \ref{LT_extensions} (iv). Then there is an embedding of $\bcO$-algebras $\cO_{\brF_n} \hra R_n$. 

\vskip5pt

{\rm (vi)} Let $\vpi_n \in \cO_{\brF_n}$ be a uniformizer. One has
\[\vpi_n \sim \prod_{[\alpha] \in \cP_n} \lambda_n^\univ(\alpha) \;.\]

\vskip5pt

{\rm (vii)}  For $\beta \in (\pi^{-n}\cO/\cO)^h_\prim$ we have the following prime factorization of $\lambda_n^\univ(\alpha)$:
\[\lambda_n^\univ(\beta) \sim \prod_{[\alpha] \in \cP_{n+1}, [\pi\alpha] = [\beta]} \lambda_{n+1}^\univ(\alpha)^q \;.\]

\vskip5pt

{\rm (viii)} For every $[v] \in \bbP^{h-1}(\cO)$ the sequence of ideals $\Big((\lambda_n^\univ(\pi^{-n}v + \cO^h))\Big)_{n>0}$ satisfies 
the conditions in \ref{prime_ideals}.

\vskip5pt

{\rm (ix)} Conversely, for every sequence of prime ideals $(\frp_n)_{n>0}$ as in \ref{prime_ideals} there is a unique $[v] \in \bbP^{h-1}(\cO)$ such that $(\lambda_n^\univ(\pi^{-n}v + \cO^h)) = \frp_n$ for all $n>0$.
\end{prop}

\Pf (i) The first statement is contained in \cite[4.3]{Drinfeld}. For the second statement let $S \subseteq R_n$ be the subring generated over $R_0$ by $\lambda_n^\univ(\alpha_1), \ldots, \lambda_n^\univ(\alpha_h)$. As $R_n$ is finite over $R_0$ it follows that $R_n/\frm_{R_0}R_n$ is a finite-dimensional $\kb_F$-vector space generated by finitely many monomials in the $\lambda_n^\univ(\alpha_i)$. In particular, $R_n = S + \frm_{R_0}R_n$. Hence $R_n = S$ by a corollary of Nakayama's Lemma.

\vskip8pt

(ii) Follows easily from (i), cf. \cite[4.2 (i)]{StrauchGeometricallyConnectedComponents}.

\vskip8pt

(iii) We fix a coordinate $T$ on $X^\univ$. Then the multiplication by $a \in \cO$ on $X^\univ$ is given by a power series $[a]_{X^\univ}(T) = aT + T^2P(T)$ with a power series $P(T) \in R_0\llbracket T \rrbracket$. If $a$ is a unit in $
\cO$, then we see that $[a]_{X^\univ}(x)/x = a + xP(x)$ is a unit in $R_n$ for all non-zero $x \in \frm_{R_n}$. It follows that 
$$\lambda_n^\univ(a\alpha) = [a]_{X^\univ}\Big(\lambda_n^\univ(\alpha)\Big) = \lambda_n^\univ(\alpha) \cdot \frac{[a]_{X^\univ}\Big(\lambda_n^\univ(\alpha)\Big)}{\lambda_n^\univ(\alpha)}$$

\vskip8pt

is associate to $\lambda_n^\univ(\alpha)$. This shows that $\lambda_n^\univ(\alpha) \sim \lambda_n^\univ(\beta)$ if $[\alpha] = [\beta]$ in $\cP_n$. The converse is in \cite[4.2 (i)]{StrauchGeometricallyConnectedComponents}.

\vskip8pt

(iv) Suppose $[\pi \alpha] = [\beta]$. By (iii) we have $\lambda_n^\univ(\beta) \sim \lambda_n^\univ(\pi \alpha)$, and 
$$\lambda_n^\univ(\pi \alpha) = \lambda_{n+1}^\univ(\pi \alpha) = [\pi]_{X^\univ}\Big(\lambda_{n+1}^\univ(\alpha)\Big) = \lambda_{n+1}^\univ(\alpha) \cdot  \frac{[\pi]_{X^\univ}\Big(\lambda_{n+1}^\univ(\alpha)\Big)}{\lambda_{n+1}^\univ(\alpha)} \;,$$

\vskip8pt

hence $\lambda_{n+1}^\univ(\alpha)$ divides $\lambda_n^\univ(\beta)$ in $R_{n+1}$. In particular, $(\lambda_{n+1}^\univ(\alpha)) \cap R_n = (\lambda_n^\univ(\pi \alpha))$.
\vskip8pt

Now suppose $\lambda_{n+1}^\univ(\alpha)$ divides $\lambda_n^\univ(\beta)$ in $R_{n+1}$. Then $(\lambda_{n+1}^\univ(\alpha)) \cap R_n = (\lambda_n^\univ(\beta))$. On the other hand, we have just seen shown that $(\lambda_{n+1}^\univ(\alpha)) \cap R_n = (\lambda_n^\univ(\pi \alpha))$. Therefore, $\lambda_n^\univ(\beta) \sim \lambda_n^\univ(\pi \alpha)$. Now we use (iii) to conclude. 

\vskip8pt

(v) and (vi) These assertions are in \cite[3.4, 4.2 (ii)]{StrauchGeometricallyConnectedComponents}. 

\vskip8pt

(vii) We apply statement (vi) twice, for $n$ and for $n+1$, and obtain:

$$\prod_{[\beta] \in \cP_n} \lambda_n^\univ(\beta) \sim \vpi_n \sim \vpi_{n+1}^q \sim \prod_{[\alpha] \in \cP_{n+1}} \lambda_{n+1}^\univ(\alpha)^q = \prod_{[\beta] \in \cP_n} \prod_{[\alpha] \in \cP_{n+1}, [\pi\alpha] = [\beta]} \lambda_{n+1}^\univ(\alpha)^q\;.$$

\vskip8pt

Assertion (vii) follows now from statement (iv). 

\vskip8pt

(viii) In $\brF_1$ we have $\pi \sim \vpi_1^{q-1}$, and thus $\pi \sim  \prod_{[\alpha] \in \cP_1} \lambda_1^\univ(\alpha)^{q-1}$, by (vi). This shows the first condition in \ref{prime_ideals}. The second condition now follows from statement (iv). 

\vskip8pt

(ix) As we have seen in the proof of (viii), any principal prime ideal $\frp_1$ of $R_1$ dividing $(\pi)$ must be generated by one of $\lambda_1(\alpha_1)$ for a unique $[\alpha_1] \in \cP_1$. By (vii), any principal prime ideal $\frp_{n+1}$ of $R_{n+1}$ dividing $(\lambda_n(\alpha_n))$, with $\alpha_n \in \cP_n$, must be generated by an element $\lambda_{n+1}(\alpha_{n+1})$ with $\alpha_{n+1} \in \cP_{n+1}$ and $[\pi \alpha_{n+1}] = [\alpha_n]$.  One can choose elements $\tilde{\alpha}_n \in \pi^{-n}\cO^h$ such that $\tilde{\alpha}_n +\cO^h = \alpha_n$ and $\pi \tilde{\alpha}_{n+1} + \cO^h = \alpha_n$. It is easily seen that the limit $v = \lim_{n \ra \infty} \pi^n \tilde{\alpha}_n$ exists and is an element in $\cO^h_\prim$, and $\pi^{-n}v + \cO^h = \alpha_n$ for all $n>0$. This proves statement (ix). \qed

\vskip8pt

\begin{cor} The prime ideals of height one of $R_\infty = \bigcup_n R_n$ lying over $(\pi)$ are naturally parametrized by elements in $\bbP^{h-1}(\cO)$, and the action of $\GL_h(\cO)$ on the set of those prime ideals of $R_\infty$ is transitive.  
\end{cor}

\Pf This is an immediate consequence of \ref{basic} (viii) and (ix). \qed

\vskip8pt

\begin{convention} In the remainder of this section we will describe certain Galois groups. Their description will involve terms like $1+ \pi^m\cO/(\pi^n)$ or $1+\pi^mM_{h-1}(\cO/(\pi^n))$, for $n \ge m \ge 0$. When $m=0$ we will interpret these terms as meaning $(\cO/(\pi^n))^\times$ and $\GL_{h-1}(\cO/(\pi^n))$, respectively. 
\end{convention}

\vskip8pt

\begin{para}{\it The fields $\cK_n$.} We denote by $\cK_n$ and $\cK_\infty$ the fields of fractions of $R_n$ and $R_\infty$, respectively. Furthermore, we let $e_1 = (1, 0, \ldots, 0), \ldots, e_h = (0, \ldots, 0, 1)$ be the standard generators of $\cO^h$. 
\end{para}

\vskip8pt

\begin{cor}\label{galois1} {\rm (i)} For every $\sigma \in \Gal(\cK_n \midc \cK_0)$ there is a unique matrix $(a_{i,j})_{1 \le i,j, \le h} \in \GL_h(\cO/(\pi^n))$ such that for all $j = 1, \ldots, h$: 

\vskip8pt

$\sigma\Big(\lambda^\univ_n(\pi^{-n}e_j + \cO^h)\Big) = [a_{1,j}]_{X^\univ}\Big(\lambda^\univ_n(\pi^{-n}e_1 + \cO^h)\Big) +_{X^\univ} \ldots$

\hfill $ \ldots +_{X^\univ} \; [a_{h,j}]_{X^\univ}\Big(\lambda^\univ_n(\pi^{-n}e_h + \cO^h)\Big)$

\vskip8pt

The so-defined map $\gamma: \Gal(\cK_n \midc \cK_0) \ra \GL_h(\cO/(\pi^n))$ is an isomorphism. 

\vskip8pt

{\rm (ii)} For varying $n \ge m \ge 0$, the isomorphism $\gamma$ in {\rm (i)} is compatible with the obvious transition maps 

$$\Gal(\cK_n \midc \cK_0) \lra \Gal(\cK_m \midc \cK_0) \;,$$
$$\GL_h(\cO/(\pi^n)) \lra \GL_h(\cO/(\pi^m)) \;,$$

\vskip8pt

and thus induces an isomorphism $\Gal(\cK_\infty \midc \cK_0) \stackrel{\simeq}{\lra} \GL_h(\cO)$.

\vskip8pt

{\rm (iii)} For $n \ge m \ge 0$, the isomorphism $\gamma$ in {\rm (i)} induces an isomorphism  

$$\Gal(\cK_n \midc \cK_m) \stackrel{\simeq}{\lra}  1+ \pi^m M_h(\cO/(\pi^n)) \;.$$ 

\vskip8pt

\end{cor}

\Pf The first statement follows directly from \ref{thm:universal-deformation} (iv) and \ref{basic} (i). The remaining statements are an easy consequence of the first. 
\qed

\vskip8pt

\begin{para}\label{defn_K'_n}{\it The fields $K'_n$.} In the following it will be convenient to fix a particular sequence of prime $\frp_n$ satisfying \ref{prime_ideals}, namely the sequence corresponding to the first standard basis vector $v = e_1 = (1, 0 ,\ldots, 0)$ by \ref{basic} (viii), i.e., $\frp_n = (\pi_n)$ where $\pi_n = \lambda^\univ_n(\pi^{-n}e_1 + \cO^h)$. Let $R'_n$ be the completion of the localization $(R_n)_{\frp_n}$ (with respect to the topology defined by the maximal ideal), denote by $K'_n = \Frac(R'_n)$ its field of fractions, and by $k'_n$ the residue field of $K'_n$. Let $k'_{n,\sep} \subseteq  k'_n$ be the separable closure of $k'_0$ in $k'_n$. Set $K'_\infty = \bigcup_{n \ge 0} K'_n$.
\end{para}

\vskip8pt

\begin{cor}\label{cor:ramification-index-k}
The ramification index of the extension $K'_n \midc K'_0$ is $e(K'_n \midc K'_0) = (q-1) q^{n-1}$.
\end{cor}

\Pf In the proof of \ref{basic} (viii) we have shown $\pi \sim  \prod_{[\alpha] \in \cP_1} \lambda_1^\univ(\alpha)^{q-1}$ which shows that $e(K'_1 \mid K'_0) = q-1$. For $n>0$ it follows from \ref{basic} (vii) that $e(K'_{n+1} \midc K'_0) = q$. 
\qed

\vskip8pt

\begin{prop}\label{prop:galois-group-k} {\rm (i)} $R'_n$ is generated over $R'_0$ by $\lambda^\univ_n(\pi^{-n}e_1 + \cO^h), \ldots, \lambda^\univ_n(\pi^{-n}e_h + \cO^h)$, and the residue field $k'_n$ of $K'_n$ is generated over $k'_0$ by the images of $\lambda^\univ_n(\pi^{-n}e_2 + \cO^h), \ldots, \lambda^\univ_n(\pi^{-n}e_h + \cO^h)$.

\vskip8pt

{\rm (ii)} The isomorphism $\gamma: \Gal(\cK_n \midc \cK_0) \car \GL_h(\cO/(\pi^n))$ in \ref{galois1} induces an isomorphism

\[\Gal(K'_n \midc K'_0) \car
	\left\{
		\left(
		\begin{array}{c|ccc}
		a_{1,1} & a_{1,2} & \cdots & a_{1,h}\\
		\hline
		0\\
		\vdots & & A'\\
		0
		\end{array}
		\right)\ 
	\middle| \ 
		{\renewcommand*{\arraystretch}{1.4}
		\begin{array}{l}
		a_{1,1}\in(\cO/\pi^n \cO)^\times,\\
		a_{1,j}\in \cO/(\pi^n) \; {\rm for} \; j>1,\\
		A'\in \GL_{h-1}(\cO/(\pi^n))
		\end{array}
		}
	\right\} \;,\]
	
\vskip8pt

which we again denote by $\gamma$.

\vskip8pt

{\rm (iii)} Let $I(K'_n \midc K'_0) \subseteq \Gal(K'_n \midc K'_0)$ be the inertia subgroup. Then the isomorphism $\gamma$ in {\rm (ii)} maps $I(K'_n \midc K'_0)$ isomorphically onto the subgroup  

\[\left\{
		\left(
		\begin{array}{c|ccc}
		a_{1,1} & a_{1,2} & \cdots & a_{1,h}\\
		\hline
		0\\
		\vdots & & I_{h-1}\\
		0
		\end{array}
		\right)\ 
	\middle| \ 
		{\renewcommand*{\arraystretch}{1.4}
		\begin{array}{l}
		a_{1,1}\in(\cO/\pi^n \cO)^\times,\\
		a_{1,j}\in \cO/(\pi^n) \; {\rm for} \; j>1\\
		\end{array}
		}
	\right\} \;,\]
	
\vskip8pt

where $I_{h-1}$ is the $(h-1)\times (h-1)$-identity matrix. Furthermore, the residue field extension $k'_n \midc k'_0$ is normal and the subextsion $k'_{n,\sep} \midc k'_0$ is Galois. The isomorphism $\gamma$ in {\rm (ii)} induces an isomorphism $\Aut(k'_n \midc k'_0) = \Gal(k'_{n,\sep} \midc k'_0) \car \GL_{h-1}(\cO/(\pi^n))$.

\vskip8pt

{\rm (iv)} For any $n \ge m \ge 0$ the isomorphism $\gamma$ in {\rm (ii)} induces an isomorphism

\[\Gal(K'_n \midc K'_m) \simeq
	\left\{
		\left(
		\begin{array}{c|ccc}
		a_{1,1} & a_{1,2} & \cdots & a_{1,h}\\
		\hline
		0\\
		\vdots & & A'\\
		0
		\end{array}
		\right)\ 
	\middle| \ 
		{\renewcommand*{\arraystretch}{1.4}
		\begin{array}{l}
		a_{1,1}\in(1+\pi^m\cO)/(1+\pi^n \cO),\\
		a_{1,j} \in (\pi^m)/(\pi^n) \; {\rm for} \; j>1,\\
		A' \in I_{h-1}+\pi^m M_{h-1}(\cO/(\pi^n))
		\end{array}
		}
	\right\} \;.\]

\vskip8pt

{\rm (v)} The isomorphism $\gamma$ in {\rm (ii)} induces, for every $m \in \bbN$, an isomorphism

\[\Gal(K'_\infty \midc K'_m) \simeq
	\left\{
		\left(
		\begin{array}{c|ccc}
		a_{1,1} & a_{1,2} & \cdots & a_{1,h}\\
		\hline
		0\\
		\vdots & & A'\\
		0
		\end{array}
		\right)\ 
	\middle| \ 
		{\renewcommand*{\arraystretch}{1.4}
		\begin{array}{l}
		a_{1,1}\in 1+ \pi^m \cO,\\
		a_{1,j}\in \pi^m\cO \; {\rm for} \; j>1,\\
		A'\in I_{h-1}+\pi^m M_{h-1}(\cO)
		\end{array}
		}
	\right\} \;.\]
	
\vskip8pt
\end{prop}

\Pf (i) Follows from \ref{basic} (i). 

\vskip8pt

(ii) This is statement \cite[4.1 (ii)]{StrauchGaloisOnTorsion}. The integer $n$ (resp. $m$, resp. $h$) in this reference corresponds to $h$ (resp. $n$, resp. $1$) here. The universal Drinfeld level structure is denoted by $\phi$ in loc.\ cit., and the ideal $\frp_{1,m}$ of loc.\ cit.\ corresponds to the ideal $(\pi_n)$ considered here. 

\vskip8pt

(iii) This is statement \cite[4.1 (iii)]{StrauchGaloisOnTorsion}. It is a general fact that the residue field extension $k'_n \midc k'_0$ is normal (and the extension $k'_{n,\sep} \midc k'_0$ is therefore Galois), and that the map from $\Gal(K'_n \midc K'_0)/I(K'_n \midc K'_0)$ to the automorphism group $\Aut(k'_n \midc k'_0)$ is an isomorphism, cf. \cite[Theorem 2 of Ch.~5, \textsection 2.2]{Bourbaki}.

\vskip8pt

(iv) Follows straightforwardly from the isomorphism in (ii).

\vskip8pt

(v) It follows from \ref{galois1} that, for fixed $m$, the isomorphism in (iv) is compatible with the natural transition maps on both sides, as $n$ varies. Passing to the projective limit shows the assertion.
\qed

\vskip8pt

\begin{para}\label{defining_Kn}{\it The fields $\tK_n$ and $K_n$.} For $\ell \in \bbN$ let $K'_{\ell,u} \subseteq K'_\ell$ be the maximal subextension of $K'_\ell \midc K'_0$ which is unramified over $K'_0$. The residue field of $K'_{\ell,u}$ is the separable closure $k'_{\ell,\sep}$ of $k'_0$ in $k'_\ell$. Set $\tK_0 = \bigcup_{\ell \ge 0} K'_{\ell,u}$, and define $K_0$ to be the $\pi$-adic completion of $\tK_0$. For $n \geq 0$, put 

$$\tK_n = K'_n \tK_0 \;\; \mbox{ and } \;\; K_n = K'_n K_0 \;.$$

\vskip8pt

The fields $\tK_n$ and $K_n$ are discretely valued and $K_n$ is complete. Denote by $A_n$ (resp. $\tA_n$) the ring of integers of $K_n$ (resp. $\tK_n$) and by $k_n$ (resp. $\tk_n$) its residue field. As completion does not affect the residue field, the canonical map $\tk_n \ra k_n$ is an isomorphism.
\end{para}

\vskip8pt

\begin{rem}\label{GaloisNonGalois} For $n \ge m \ge 0$ the extension $\tK_n \midc \tK_m$ is Galois of finite degree, and $\tK_n$ is also a Galois extension of any of the fields $K'_{\ell,u}$ (this extension is of infinite degree). Similarly, the extension $K_n \midc K_m$ is Galois of finite degree. But, if $h>1$, the $p$-adically complete field $K_n$ is neither Galois over $K'_m$ ($m \le n$), nor is it Galois over $\tK_m$ ($m \le n$) or $K'_{\ell,u}$ (any $\ell$).
\end{rem} 

\vskip8pt

\begin{prop}\label{prop:gal-groups-knprime} 
{\rm (i)} For $\ell' \ge n \ge \ell \ge m \ge 0$, the isomorphism $\gamma$ in \ref{prop:galois-group-k} (ii) induces an isomorphism

\[\Gal(K'_n K'_{\ell',u} \midc K'_m K'_{\ell,u}) \simeq
	\left\{
		\left(
		\begin{array}{c|ccc}
		a_{1,1} & a_{1,2} & \cdots & a_{1,h}\\
		\hline
		0\\
		\vdots & & A'\\
		0
		\end{array}
		\right)\ 
	\middle| \ 
		{\renewcommand*{\arraystretch}{1.4}
		\begin{array}{l}
		a_{1,1}\in(1+\pi^m\cO)/(1+\pi^n \cO),\\
		a_{1,j} \in (\pi^m)/(\pi^n) \; {\rm for} \; j>1,\\
		A' \in I_{h-1}+\pi^\ell M_{h-1}(\cO/(\pi^{\ell'}))
		\end{array}
		}
	\right\} \;,\]
	
\vskip8pt

which we again denote by $\gamma$. This isomorphism is compatible with the natural transition maps on both sides as $\ell' \ge n \ge \ell \ge m$ vary.

\vskip8pt

{\rm (ii)} For any $n$ and $\ell$ the isomorphism $\gamma$ in {\rm (i)} induces an isomorphism

\[\Gal(\tK_n \midc K'_{\ell,u}) \simeq
	\left\{
		\left(
		\begin{array}{c|ccc}
		a_{1,1} & a_{1,2} & \cdots & a_{1,h}\\
		\hline
		0\\
		\vdots & & A'\\
		0
		\end{array}
		\right)\ 
	\middle| \ 
		{\renewcommand*{\arraystretch}{1.4}
		\begin{array}{l}
		a_{1,1}\in (\cO/(\pi^n))^\times,\\
		a_{1,j} \in \cO/(\pi^n) \; {\rm for} \; j>1,\\
		A' \in I_{h-1}+\pi^\ell M_{h-1}(\cO))
		\end{array}
		}
	\right\} \;.\]
	
\vskip8pt

{\rm (iii)} For any $n \ge m$ the isomorphism $\gamma$ in {\rm (i)} induces an isomorphism

$$\Gal(\tK_n \mid \tK_m) \simeq 1+ \pi^m \cO/(\pi^n) \ltimes ((\pi^m)/(\pi^n))^{h-1} \;,$$

\vskip8pt

{\rm (iv)} For all $n \ge m \ge 0$, the extension $K_n \midc K_m$ is Galois and the restriction map $\Gal(K_n \mid K_m) \ra \Gal(\tK_n \mid \tK_m)$ is an isomorphism.
We thus have $[K_n:K_m] = q^{(n-m)h}$, if $m>0$, and $[K_n:K_0] = (q-1)q^{nh-1}$. 
\vskip8pt

{\rm (v)} For all $n \ge 1$ one has $[k_n:k_{n-1}] = q^{h-1}$.
\end{prop}

\Pf (i) We consider $K'_n K'_{\ell',u}$ as a subfield of $K'_{\ell'}$. Then we have the short exact sequence

\[1 \lra \Gal(K'_{\ell'} \midc K'_n K'_{\ell',u}) \lra \Gal(K'_{\ell'} \midc K'_0) \lra \Gal(K'_n K'_{\ell',u} \midc K'_0) \lra 1 \;.\]

\vskip8pt

Consider $\sigma \in \Gal(K'_{\ell'} \midc K'_0)$ and write 

\[\gamma(\sigma) = \left(
		\begin{array}{c|ccc}
		a_{1,1} & a_{1,2} & \cdots & a_{1,h}\\
		\hline
		0\\
		\vdots & & A'\\
		0
		\end{array} \right)\]
		
\vskip8pt

with $a_{1,1} \in \cO/(\pi^{\ell'})$, $a_{1,j} \in \cO/(\pi^{\ell'})$, for $j>1$, and $A' \in \GL_{h-1}(\cO/(\pi^{\ell'})$. For $\sigma$ to act trivially on $K'_n K'_{\ell',u}$ it must, in particular, act trivially on $K'_{\ell',u}$. By \ref{prop:galois-group-k} (iii), the matrix $A'$ must therefore be the identity matrix. Since $K'_n$ is generated by the elements $\lambda^\univ_n(\pi^{-n}e_j + \cO^h)$, for $j=1, \ldots, h$, we must also have $a_{1,1} \in 1+(\pi^n)/(\pi^{\ell'})$ and $a_{1,j} \in (\pi^n)/(\pi^{\ell'})$ for $j>1$. This proves the assertion when $\ell = m = 0$. Using similar arguments we see that the subgroup $\Gal(K'_n K'_{\ell',u} \midc K'_m K'_{\ell,u})$ of $\Gal(K'_n K'_{\ell',u} \midc K'_0)$ is mapped by $\gamma$ to the group as stated. 

\vskip8pt

(ii) This follows from (i) when we take $m=0$ and when pass to the projective limit as $\ell' \ra \infty$.

\vskip8pt

(iii) This statement follows from (ii) when we take $\ell = 0$, and from the exact sequence

\[1 \lra \Gal(\tK_n \midc \tK_m) \lra \Gal(\tK_n \midc \tK_0) \lra \Gal(\tK_m \midc \tK_0) \lra 1 \;.\]

\vskip8pt

(iv) The extension $K'_n \midc K'_0$ is Galois, and so is the extension $K_n = K'_n K_0 \midc K_m = K'_m K_0$ that we obtain by taking the composite fields with $K_0$. The restriction map $\Gal(K_n \midc K_m) \ra  \Gal(\tK_n \midc \tK_m)$ is injective, because if $\sigma \in \Gal(K_n \midc K_m)$ acts trivially on $\tK_n$, which is dense in $K_n$ (for the $p$-adic topology), then it acts trivially on $K_n$. But this map is also surjective because any Galois automorphism in $\Gal(\tK_n \midc \tK_m)$ extends continuously to an automorphism of $K_n$ over $K_m$, because $\tK_n$ is dense in $K_n$.

\vskip8pt

(v) By (iii) and (iv) we have $[K_n:K_0] = (q-1)q^{n(h-1)+n-1}$ for all $n \ge 1$. For the ramification index we have $e(K_n \midc K_0) = e(\tK_n \midc \tK_0) = e(K'_n \midc K'_0) = (q-1)q^{n-1}$; cf. \ref{cor:ramification-index-k}, since $K_0$ is the completion of the unramified extension $\tK_0$ of $K'_0$. It follows that $[k_n:k_0] = q^{n(h-1)}$, from which we conclude that $[k_n:k_{n-1}] = q^{h-1}$. \qed

\section{Strictly deeply ramified towers}\label{deep_ramification}\setcounter{subsection}{1}

In the rest of this paper we will only consider the tower $K_\bullet$ constructed in section \ref{tower} when $\cO = \Zp$. In particular, we have 

\[q=p \;\; \mbox{and} \;\; \pi = p \;.\]

In the following we will always write $p$ instead of $q$, when using formulas from the preceding section (involving cardinalities), but we keep writing $\pi$ instead of $p$ when referring to the uniformizer of $\cO = \Zp$.

\vskip8pt

The reason for restricting our attention to the case of $\Zp$ is because of the way the theory of strictly deeply ramified towers of fields has been developed by Scholl. We note, however, that Scholl's theory can be generalized to a setting which would allow us to work here with the ring of integers $\cO$ of a finite extension of $\Qp$, cf. \cite[2.3.1]{Scholl}. 

\vskip8pt

For $i=1,\ldots,h$ we put $Y_{n,i} = \lambda^{\rm univ}_n(p^{-n}e_i + \bbZ_p^h)$ which we consider as elements of $A_n$. Recall that, with this definition,  $\pi_n = Y_{n,1}$, cf. \ref{defn_K'_n}. 

\vskip8pt

Recall that we denote by $k_n$ the residue field of $K_n = K'_n K_0$, cf. \ref{defining_Kn}, and by $k'_n$ the residue field of $K'_n$, cf. \ref{defn_K'_n}. By \ref{basic} the elements $Y_{n,1}, \ldots, Y_{n,h}$ form a system of parameters of $R_n$, and generate $R_n$ as $R_0$-algebra. Therefore, $k'_n = \Frac(R_n/\pi_n R_n)$ is thus generated by 

\[y_{n,i} := Y_{n,i} \mbox{ mod } \frp_n \;, i=2, \ldots, h \;,\]

\vskip8pt

over $k'_0$. 

\begin{prop}\label{minpol}
For all $n \ge 1$ and $i=1, \ldots, h$, the minimal polynomial of $Y_{n,i}$ over $K_{n-1}$ is
    \[ Q_i(T) = Q_{n,i}(T) = \prod_{a\in\Fp} (T - (Y_{n,i} +_{X^\univ} [a]_{X^\univ}(Y_{1,1}))) \;,\]

and has coefficients in $A_{n-1}$, except if $n=1=i$ in which case 

\[Q_{1,1}(T) = \prod_{a\in\Fp, \,a \neq -1} (T - (Y_{1,1} +_{X^\univ} [a]_{X^\univ}(Y_{1,1}))) \;.\]

For $2 \le i \le h$, the reduction of $Q_{n,i}$ modulo the maximal ideal $(\pi_{n-1})$ of $A_{n-1}$ is $T^p-y_{n,i}^p$. 

\end{prop}
\Pf
We only treat the case $n>1$; the case $n=1$ is very similar (with the obvious modifications). By \ref{prop:gal-groups-knprime} (iii, iv), the Galois group of $K_n \mid K_{n-1}$ is isomorphic to the group \[(1+\pi^{n-1}\cO)/(1+\pi^n \cO)\ltimes (\pi^{n-1}\cO/\pi^n \cO)^{h-1} \;,\] 

and the Galois action is given by the formula in \ref{galois1}. This means the Galois conjugates of $Y_{n,i}$ are elements \[Y_{n,i} +_{X^\univ} [a\pi^{n-1}]_{X^\univ}(Y_{n,1})\] for $a\in\Fp$. This shows that the minimal polynomial of $Y_{n,i}$ over $K_{n-1}$ is thus
    \[ Q_i(T) = \prod_{a\in\Fp} (T - (Y_{n,i} +_{X^\univ} [a\pi^{n-1}]_{X^\univ}(Y_{n,1}))) \;. \] 
    
Since all roots are in $A_n$, the coefficients of this polynomial are in $A_n \cap K_{n-1} = A_{n-1}$. 

\vskip8pt

The universal formal group law, as any one-dimensional formal group law, has the property that 
\[T_1 +_{X^\univ} T_2 = T_1 + T_2 + T_1T_2 \cdot (\cdots)\]

This implies that $Q_{n,i}(T) \equiv \prod_{a \in \Fp} (T-(Y_{n,i} \mbox{ mod } \pi_n) = T^p-y_{n,i}^p \mbox{ mod } \pi_n$. As the coefficients of $Q_{n,i}$ are in $A_{n-1}$, we also have $Q_{n,i}(T) \mbox{ mod } \pi_{n-1} = T^p - y_{n,i}^p$.
\qed

\begin{prop}\label{generators}
\begin{enumerate}
    \item The residue field $k_n$ is generated as a field over $k_{n-1}$ by the elements $y_{n,2},\ldots,y_{n,h}$,, and $\{y_{n,2}^{i_2} \cdots  y_{n,h}^{i_h} \midc 0 \le i_j \le p-1\}$ is a basis of $k_n$ over $k_{n-1}$.
    \item The elements $Y_{n,1},\ldots,Y_{n,h}$ generate $A_n$ as an algebra over $A_{n-1}$, and $A_n$ is a free $A_{n-1}$-module with basis \[\{Y_{n,1}^{i_1}\cdots Y_{n,h}^{i_h} \midc \forall j \in \{1, \ldots,h\}:  0\leq i_j\leq p-1\} \;,\] 
    except if $n=1$ in which case $A_1$ is free over $A_0$ with basis $\{Y_{n,1}^{i_1}\cdots Y_{n,h}^{i_h} \midc \forall j \in \{1, \ldots,h\}:  0 \le i_1 \le p-2, \, 0 \leq i_j\leq p-1 \mbox{ for } j>1\}$.
\end{enumerate}
\end{prop} 
\Pf
(i) Recall that $K_0$ is the $p$-adic completion of $\tK_0$. Since $K'_n$ is finite over $K'_0$, it follows that $K_n$ is also the $p$-adic completion of $\tK_n = K'_n \tK_0$. By definition, $\tK_0 = \bigcup_{\ell \ge 0} K'_{\ell,u}$, cf. \ref{defining_Kn}, and we thus have $\tK_n = \bigcup_{\ell \ge 0} K'_n K'_{\ell,u}$. Because the residue field does not change after passing to the completion, the residue field $k_n$ of $K_n$ is equal to the residue field of $\tK_n$, and the residue field of $\tK_n$ is the union of the residue fields of the $K'_n K'_{\ell,u}$. By definition, $K'_{\ell,u}$ is unramified over $K'_0$ and it's residue field is the separable closure $k'_{\ell,\sep}$ of $k'_0$ in $k'_\ell$. By \cite[2.4.8]{FriedJarden} the residue field of the composite field $K'_n K'_{\ell,u}$ is thus equal to $k'_n k'_{\ell,\sep}$, the composition of the residue fields. The union of these fields is then $k'_n k_0$. By the remark before \ref{minpol}, the field $k_n'$ is generated over $k'_{n-1}$ by $y_{n,i}$, $2 \le i \le h$, and $k_n$ is thus generated over $k_{n-1}$ by those same elements. By \ref{minpol}, these elements are of degree $\le p$, and because of $[k_n:k_{n-1}] = p^{h-1}$, they are indeed of degree $p$.

\vskip8pt

(ii) We observe that $A_{n-1}$ is a local ring with maximal ideal $\pi_{n-1} A_{n-1}$. Since $\pi_{n-1}A_n = \pi_n^pA_n = Y_{n,1}^pA_n$ (except if $n=1$ in which case $\pi_0A_1 = pA_1 = \pi_1^{p-1}A_1$), we have a filtration 

\[A_n/\pi_{n-1} A_n  = A_n/(\pi_n^p) \subseteq (\pi_n^{p-1}) /(\pi_n^p) \subseteq \ldots \subseteq (\pi_n)/(\pi_n^p) \subseteq  A_n/(\pi_n^p)\] 

(and similarly when $n=1$, when we replace $\pi_n^p$ by $\pi_1^{p-1}$. The successive quotients are $(\pi_n^i)/(\pi_n^{i+1}) \simeq k_n$, hence are finitely generated as $k_{n-1}$-vector space by monomials in $y_{n,2},\ldots, y_{n,h}$. We thus see that as a $k_{n-1}$-algebra, we have 
\[A_n/\pi_{n-1} A_n = k_{n-1} [\bar{Y}_{n,1}, \ldots, \bar{Y}_{n,h}],\] 

where $\bar{Y}_{n,i}$ denotes the image of $Y_{n,i} \in A_n$ modulo $\pi_{n-1}$. By \ref{cor:ramification-index-k}, $e(K'_n\mid K'_{n-1})=p$ if $n>1$ (and is equal to $p-1$ if $n=1$), we have that $Y_{n,1}$ is of degree $p$ (resp. $p-1$, if $n=1$) over $A_{n-1}$, and each element $\bar{Y}_{n,i}$ is of degree $p$ over $k'_{n-1}$ for $i\geq 2$ as above, the $k'_{n-1}$-vector space $A'_n/\pi_{n-1} A'_n$ is generated by 
\[\{\bar{Y}_{n,1}^{i_1}\cdots \bar{Y}_{n,h}^{i_h} \midc \forall j \in \{1, \ldots,h\}:  0\leq i_j\leq p-1\} \;\] 

(if $n-1$ then it suffices that $i_1 \le p-2$). By Nakayama's Lemma, this generating set can be lifted to a generating set of $A_n$ as an $A_{n-1}$-module, which is to say that $A_n$ is generated by the set 
\[\{Y_{n,1}^{i_1}\cdots Y_{n,h}^{i_h} \midc \forall j \in \{1, \ldots,h\}:  0\leq i_j\leq p-1\}\] 

(if $n-1$ then it suffices that $i_1 \le p-2$) as an $A_{n-1}$-module. It is a general fact that $A_n$ is free over $A_{n-1}$ of degree $[K_n:K_{n-1}] = p^h$ (resp. $[K_1:K_0] = (p-1)p^{h-1}$), cf. \cite[ch. II, \S 2, Prop. 3]{Serre_LF}, and those elements must then be basis of $A_n$ as $A_{n-1}$-module.  
\qed

\begin{prop}\label{big_field} For every $n \ge 0$ we have $[k_n:k_n^p] = p^{h-1}$. Therefore, each field $K_n$ is a d-big local field in the sense of Scholl \cite[1.1]{Scholl} with $d=h-1$.
\end{prop}

\Pf Recall that $k'_0 = \overline{\bbF}_p \llbracket u_1, \ldots, u_{h-1} \rrbracket$. Thus the extension $k'_0 \mid (k'_0)^p$ has degree $p^{h-1}$, because it is generated by the elements $u_1, \ldots, u_{h-1}$, each of which is of degree $p$ over $(k'_0)^p$ and thus inseparable over $(k'_0)^p$. Via the isomorphism $k'_{n, \sep} \to (k'_{n, \sep})^p$, $x \mapsto x^p$, we have $\Gal(k'_{n, \sep} \mid k'_0) \simeq \Gal((k'_{n,\sep})^p \mid (k'_0)^p)$ for each $n$, and thus each extension $k'_{n, \sep} \mid (k'_{n, \sep})^p$ also has degree $p^{h-1}$. As the extensions $k'_{n,\sep} \mid k'_0$ are separable, so are the extensions $(k'_{n,\sep})^p \mid (k'_0)^p$, so the elements $u_1, \ldots, u_{h-1}$ do not belong to any of the fields $(k'_{n,\sep})^p$. Thus these elements generate each of the extensions $k'_{n,\sep} \mid (k'_{n,\sep})^p$ and hence the extension $k_0 \mid (k_0)^p$, where $k_0 = \bigcup_n k'_{n,\sep}$ is the residue field of $K_0$. The fact that the elements $u_1^{i_1} \cdots u_{h-1}^{i_{h-1}}$, with the exponents running independently over $0, \ldots, p-1$, are linearly independent over each field $(k'_{n,\sep})^p$ implies that they are still linearly independent over the union $(k_0)^p$, and thus the extension $k_0 \mid k_0^p$ has degree $p^{h-1}$.

\vskip8pt

Suppose towards induction that $k_{n-1} \mid k_{n-1}^p$ has degree $p^{h-1}$. By \ref{prop:gal-groups-knprime} (v), the extension $k_n \mid k_{n-1}$ has degree $p^{h-1}$, and via the isomorphism $x \mapsto x^p$, we conclude that the extension $k_n^p \mid k_{n-1}^p$ also has degree $p^{h-1}$. By \ref{generators} we have 
\[k_n^p = (k_{n-1}(y_{n,2}, \ldots, y_{n,h}))^p = k_{n-1}^p(y_{n,2}^p, \ldots, y_{n,h}^p) \;.\] 

By \ref{minpol}, the reduction of the minimal polynomial of $Y_{n,i}$ over $K_{n-1}$ to $k_{n-1}$ is equal to $T^p - y_{n,i}^p$, and thus $y_{n,i}^p \in k_{n-1}$ for each $i$. It follows that $k_n^p = k_{n-1}^p(y_{n,2}^p, \ldots, y_{n,h}^p) \subseteq k_{n-1}$. As both have the same degree over $k_{n-1}^p$, we must have $k_n^p = k_{n-1}$. It then follows that $[k_n : k_n^p] = [k_n : k_{n-1}] = p^{h-1}$. \qed

\vskip8pt

We recall the definition of a strictly deeply ramified tower.

\begin{dfn} {\rm \cite[1.3]{Scholl}} Let $d$ be a non-negative integer, and let 
\[L_\bullet = (L_0 \subseteq L_1 \subseteq L_2 \subseteq \ldots \;) \] 

be a tower of $d$-big local fields. The tower $L_\bullet$ is called {\it strictly deeply ramified} if there exists an integer $n_0 \ge 0$ and an ideal $\xi \subseteq \cO_{L_{n_0}}$ with $0 < v_p(\xi) \le 1$ such that the following condition holds: for every $n \ge n_0$ the extension $L_n/L_{n-1}$ has degree $p^{d+1}$, and there exists a surjection 
\[ \Omega_{\cO_{L_n}/\cO_{L_{n-1}}} \twoheadrightarrow (\cO_{L_n}/\xi\cO_{L_n})^{d+1} \;.\]

\end{dfn}

\vskip8pt

We now arrive at our first goal, namely the proof of Result~\ref{result1}.
\begin{prop}
\label{strictly-deeply-ramified}
The tower $(K_n)_n$ is strictly deeply ramified.
\end{prop}

\Pf
By \ref{prop:gal-groups-knprime}, we have $[K_n:K_{n-1}] = p^h$ for $n \ge 2$. It remains to show that for all $n \geq 2$, there exists a surjection $ \Omega_{A_n \midc A_{n-1}} \to (A_n/\pi A_n)^h$. By \ref{minpol}, the minimal polynomial of $Y_{n,i}$ over $K_{n-1}$ is
\begin{align*}
Q_{n,i}(T) &= \prod_{a\in\Fp} (T - (Y_{n,i} +_{X^\univ} [a\pi^{n-1}]_{X^\univ}(Y_{n,1}))) \;,
\end{align*}

which has coefficients in $A_{n-1}$. It then follows that
\begin{align*}
	\left[\frac{d}{dT}Q_{n,i}\right](Y_{n,i}) &= \sum_{a \in \Fp} \prod_{b\neq a} (Y_{n,i} - (Y_{n,i} +_{X^\univ}  \; [b]_{X^\univ}(Y_{1,1})))\\
	  &= \prod_{b\neq0} (Y_{n,i} - (Y_{n,i} +_{X^\univ}  \; [b]_{X^\univ}(Y_{1,1})))\\
	 &= \prod_{b\neq0} (Y_{n,i} - (Y_{n,i} + b\pi_1 + b\pi_1 Y_{n,i} \cdot  (\cdots)))\\
	 &= \prod_{b\neq0} (-b\pi_1)(1 + Y_{n,i} \cdot (\cdots)),
	\end{align*}
where $(\cdots)$ is an element of $A_n$. As $-b(1 + Y_{n,i}(\cdots))$ is a unit, and $K_1 / K_0$ has ramification index $p - 1$, we have $|(d/dT) Q_{n,i}(Y_{n,i})| = |\pi_1^{p-1}| = |\pi|$. In particular,

\begin{numequation}\label{derivative} \left[\frac{d}{dT}Q_{n,i}\right](Y_{n,i}) \in \pi A_n \;.
\end{numequation}

We now show that $\Omega_{A_n \mid A_{n-1}  }$ is a free $A_n/\pi A_n$-module of rank $h$. As $Q_{n,i}(Y_{n,i})=0$, we have
\[0 = d\Big(Q_i(Y_{n,i})\Big) = \Big[\frac{d}{dT} Q_{n,i}\Big](Y_{n,i}) \cdot dY_{n,i} = \varepsilon \pi dY_{n,i} \;,\]
for some unit $\varepsilon \in A_n$. Therefore $\pi dY_{n,i} = 0$. Because the elements $Y_{n,i}$ generate $A_n$ as an algebra over $A_{n-1}$, the $dY_{n,i}$ generate $\Omega_{A_n \mid A_{n-1}}$ as a module over $A_n$. Hence we have shown that $\pi$ annihilates $\Omega_{A_n \mid A_{n-1}  }$, and $\Omega_{A_n \mid A_{n-1}}$ is thus a module over $A_n/\pi A_n$. By the definition of the polynomials $Q_{n,i}$, the map
    \[\theta: A_{n-1}[T_1, \ldots, T_h]/(Q_{n,1}(T_1), \ldots, Q_{n,h}(T_h)) \lra A_n\;, \;\; T_i \mapsto Y_{n,i} \;,\]
is well-defined, and by \ref{generators} it is surjective. By \ref{minpol} the domain of $\theta$ is a free $A_{n-1}$-module of degree $[K_n:K_{n-1}] $, and so is the target of $\theta$. Therefore, $\theta$ is an isomorphism of $A_{n-1}$-algebras.  

\vskip8pt 

Let $\widetilde{\theta}$ be the composition of $A_{n-1}[T_1, \ldots, T_h] \ra A_{n-1}[T_1, \ldots, T_h]/(Q_{n,1}(T_1), \ldots, Q_{n,h}(T_h))$ and $\theta$. By \ref{derivative} we have $\theta\Big(\frac{d}{dT_i}Q_{n,i}(T_i)\Big) \in \pi A_n$. This implies that the map

\[A_{n-1}[T_1, \ldots, T_h] \stackrel{\frac{d}{d T_i}}{\lra}  A_{n-1}[T_1, \ldots, T_h] \stackrel{\widetilde{\theta}}{\lra} A_n \lra A_n/\pi A_n\]

factors via $A_{n-1}[T_1, \ldots,T_h]/(Q_{n,1}(T_1), \ldots, Q_{n,h}(T_h)) \cong A_n$ and induces a $A_{n-1}$-linear derivation $\partial_i: A_n \ra A_n/\pi A_n$. By the universal property of $\Omega_{A_n \mid A_{n-1}}$, there is a unique $A_n$-linear map $\psi_i: \Omega_{A_n \mid A_{n-1}} \ra A_n/\pi A_n$ such that $\partial_i (a) = \psi_i(d a)$ for all $a \in A_n$. 

\vskip8pt

Suppose that in $\Omega_{A_n \mid A_{n-1}}$ we have a relation $\sum a_j dY_{n,j} = 0$ for some $a_j \in A_n/\pi A_n$. Applying $\psi_i$ to this equation we find $0 = \sum a_j \psi_i(d Y_{n,j}) = a_i$. Thus $\Omega_{A_n \mid A_{n-1}}$ is a free module over $A_n/\pi A_n$, with basis $dY_{n,1},\ldots,dY_{n,h}$. The tower $(K_n)_n$ is thus strictly deeply ramified.
\qed

\section{Is the Lubin-Tate Tower a Kummer tower?}

In this section we investigate whether the Lubin-Tate tower of fields introduced above is a Kummer tower (as recalled in the introduction). We will only consider the case when the formal group $\bbX$ has height 2.  

\vskip8pt

In this section we use the following convention. When we study group cohomology, we let a group $G$ act on a abelian group $A$ from the left: $G \times A \ra A$, $(g,a) \mapsto g \cdot a$. Furthermore, we write a 1-cocycle $c$ on $G$ with values in $A$ as $g \mapsto c_g$, i.e., $c_g$ is the value of $c$ on $g \in G$. The cocyle $c$ then satisfies $c_{gh} = g \cdot c_h + c_g$.

\subsection{Preliminaries on Galois cohomology}

We begin by recalling the following elementary result:

\begin{lemma}[{{\cite[6.2.2]{Weibel}}}]
\label{lem:H1}
Let $U = \langle u \rangle$ be a cyclic group of order $d$, and let $M$ be an abelian $U$-module. There exists a canonical isomorphism
	\begin{gather*}
	\left\{m\in M\ \middle|\ \sum_{j=1}^{d-1} u^j m = 0\right\} \Big/ \Big\{um - m\ |\ m\in M\Big\} \to H^1(U,M)\\
	\bar{m} \mapsto \left[c(m)_{u^i} = \sum_{j=0}^{i-1} u^j \cdot m\right].
	\end{gather*}
Furthermore, if the action of $U$ on $M$ is trivial and $dM = 0$, then there exists a canonical isomorphism $M\simeq H^1(U,M)$.
\end{lemma}

\begin{cor}
\label{cor:nat_act_mu}
Suppose $p>2$. Let $U = (\Z / p^n \Z)^\times$ and $M = \mu_{p^n}$. The group $U$ acts on $M$ by setting $u\cdot\zeta = \zeta^u$. Then $H^1(U,M) = 0$.
\end{cor}

\Pf
It is well known (and easy to prove) that the group $(\Z/p^n\Z)^\times$ is isomorphic to $\mu_{p-1} \times \Z/p^{n-1}\Z$ (if $p$ is odd), and it is thus itself cyclic of order $d = (p - 1) p^{n - 1}$. Let $u \in \Z \,\setminus\, p\Z$ be such that $u + p^n\Z \in U$ is a generator. Since $u \not \equiv 1 \mbox{ mod } p$, we have $p \nmid u - 1$. On the other hand, $p^n \mid u^d - 1$, so for all $\zeta \in \mu_{p^n}$, we have 

\[\prod_{j=0}^{d-1} \zeta^{u^j} = \zeta^{(u^d - 1) / (u - 1)} = 1 \;,\]

\vskip8pt

and so

\[\left\{ \zeta \in \mu_{p^n}\ \middle|\ \prod_{j = 0}^{d-1} u^j \cdot \zeta = 1\right\} = M \;.\]

\vskip8pt
	
On the other hand, as $p\nmid u-1$, we have $u - 1 \in (\Z/p^n\Z)^\times$, so

\[\{ \zeta^{u-1}\ |\ \zeta\in\mu_{p^n}\} = M \;.\]

\vskip8pt

Thus by Lemma~\ref{lem:H1}, we have $H^1(U,M) = 0$.
\qed

\begin{prop}
\label{prop:H1_semidir}
Let $\mu$ be a finite cyclic group of order $k$ which we will write multiplicatively. Put $U = (\Z / k\Z)^\times$ and $E = \Z / k\Z$, and let $G = U \ltimes E$ be the semi-direct product of $U$ with $E$, with multiplication given by

\[(\baru_1,\bare_1) \cdot (\baru_2,\bare_2) = (\baru_1 \baru_2, \baru_1 \bare_2 + \bare_1) \;,\]

\vskip8pt

where $\barx$ denotes the class of $x \in \Z$ modulo $k$. The group $G$ acts on $\mu$ via $U$, with $(\baru,\bare) \cdot \zeta = \zeta^u$ for $(\baru,\bare) \in G = (\Z/k\Z)^\times \ltimes \Z/k\Z$, $\zeta \in \mu$. Then there exists a split exact sequence \[1 \to H^1(U,\mu) \to H^1(G,\mu) \to \mu \to 1,\] where the splitting is given by mapping an element $\zeta\in\mu$ to the cohomology class of the 1-cocycle $\tilde{c}(\zeta)_{(\baru,\bare)} = \zeta^e$. Further, if $k = p^n$ for some prime $p > 2$, then $H^1(G,\mu) \simeq \mu$.
\end{prop}
\Pf
We begin with the inflation-restriction exact sequence \[0 \to H^1(U,\mu) \to H^1(G,\mu) \to H^0(U, H^1(E,\mu)) \;,\] 
cf. \cite[6.8.3]{Weibel}. As $E$ acts trivially on $\mu$, we have, by Lemma~\ref{lem:H1}, $H^1(E,\mu) \simeq \mu$, where the element $\zeta \in \mu$ corresponds to the cocycle $c(\zeta)_{\bare} = \zeta^e$, where $\bare \in E = \Z/k\Z$. The group $U$ acts on $H^1(E,\mu)$ by

\[(u \cdot c(\zeta))_{\bare} = u \cdot c(\zeta)_{\overline{u^{-1} e}} = (\zeta^{u^{-1} e})^u = \zeta^e = c(\zeta)_{\bare} \;,\]

\vskip8pt

which is to say that the action of $U$ on $H^1(E,\mu)$ is trivial. Thus 

\[H^0(U, H^1(E,\mu)) = H^1(E,\mu) \simeq \mu \;.\]

\vskip8pt

Define the splitting map as in the statement of the proposition. We check that it satisfies the cocycle condition:

\[\tilde{c}(\zeta)_{(\baru_1,\bare_1)(\baru_2,\bare_2)} = \tilde{c}(\zeta)_{(\baru_1 \baru_2, \baru_1 \bare_2 + \bare_1)} = \zeta^{u_1 e_2 + e_1} = (\tilde{c}(\zeta)_{(\baru_2,\bare_2)})^{u_1} \cdot \tilde{c}(\zeta)_{(\baru_1,\bare_1)} \;.\]

\vskip8pt

It is straightforward to check that the map $\mu \ra H^1(G,\mu)$, $\zeta \mapsto \tilde{c}(\zeta)$, is a group homomorphism and that it is a right inverse for the map $H^1(G,\mu) \to \mu$. Finally, if $k = p^n$, then by Corollary~\ref{cor:nat_act_mu} we have $H^1(U,\mu) = 0$, and thus $H^1(G,\mu) \simeq \mu$.
\qed

\begin{prop}\label{prop_Kummer}
Let $k$ be a positive integer. Suppose $L\mid K$ is a Galois extension of fields with Galois group $G = U \ltimes E$, where $U = (\Z / k\Z)^\times$ and $E = \Z / k\Z$, and the multiplication in $G$ is given by  $(\baru_1,\bare_1) \cdot (\baru_2,\bare_2) = (\baru_1 \baru_2,\baru_1 \bare_2 + \bare_1)$, where $\barx = x \mbox{ mod } k$. Suppose $L$ contains a primitive $k$-th root of unity (and therefore all $k$-th roots of unity), and suppose $G$ acts on the group $\mu_k$ of $k$-th roots of unity by $(\baru,\bare) \cdot \zeta = \zeta^u$. Then there exists a $t\in K^\times$ such that $L = K(\mu_k, t^{1/k})$.
\end{prop}
\Pf
 Let $M = L^E$. Then $\Gal(L \midc M) = E$ is cyclic of order $k$. As $E$ acts trivially on $\mu_k$, we have $\mu_k \subseteq M$. By Kummer theory, $L\mid M$ is a Kummer extension of the form $L = M(t^{1/k})$ for some $t\in M^\times$. We now want to show that one can find such an element $t$ already in $K$.

\vskip8pt

It suffices to show that there exists a $t\in K^\times$ which is an $k$-th power in $L$, $t = s^k$, on which the Galois group acts by $(\baru,\bare) \cdot s = \zeta^e s$ for some primitive $k$-th root of unity $\zeta$, since in this case $[M(s) : M] = |E|$, and thus $M(s) = L$.

\vskip8pt

Consider the exact sequence of $G$-modules 

\[1 \to \mu_k \to L^\times \to (L^\times)^k \to 1.\] 

\vskip8pt

From this we get the sequence of cohomology groups 

\[H^0(G, L^\times) \to H^0(G,(L^\times)^k) \to H^1(G,\mu_k) \to H^1(G,L^\times).\]

\vskip8pt

But the zero-th cohomology group on the left is just $K^\times$, and by Hilbert's Theorem 90 the group $H^1(G,L^\times)$ is trivial, so we obtain the sequence 

\[K^\times \to (L^\times)^k \cap K^\times \to H^1(G,\mu_k) \to 1 \;.\]

\vskip8pt

In particular, the map $(L^\times)^k \cap K^\times \to H^1(G,\mu_k)$ is surjective. As in Proposition~\ref{prop:H1_semidir}, the map $(\baru,\bare)\mapsto \zeta^e$ is a 1-cocycle, and we have a group homomorphism 

\[\mu_k \to H^1(G,\mu_k) \;,\;\; \zeta \mapsto \tilde{c}(\zeta) = [(\baru,\bare) \mapsto \zeta^e] \;.\]

\vskip8pt

Suppose that there exists a $\xi \in \mu_k$ such that $\xi^e = \xi^{u - 1}$ for all $(\baru,\bare) \in G$, i.e., the map $(\baru,\bare) \mapsto \xi^e$ is a 1-coboundary. Then $\xi^1 = \xi^{1-1} = 1$. Thus the map $\zeta \mapsto \tilde{c}(\zeta)$ is injective.

\vskip8pt

Let $\zeta$ be a primitive $k$-th root of unity. As the map $(L^\times)^k \cap K^\times \to H^1(G,\mu_k)$ is surjective, there exists a $t \in (L^\times)^k \cap K^\times$ which maps to the 1-cocycle $\tilde{c}(\zeta)$. Let $s' \in L^\times$ be such that $(s')^k = t$. By definition, under the map 

\[(L^\times)^k \cap K^\times  = H^0(G,(L^\times)^k) \to H^1(G,\mu_k) \;,\]

\vskip8pt

the element $t$ maps to the cohomology class of the 1-cocyle $\left[g \mapsto \frac{g(s')}{s'} \right]$ ($g \in G$). Thus $\tilde{c}(\zeta)$ and $\left[g \mapsto \frac{g(s')}{s'}\right]$ must be equal up to a coboundary $[g \mapsto \frac{g(\xi)}{\xi}]$, for some $\xi \in \mu_k$. That is to say: $\frac{g(s')}{s'} = \frac{g(\xi)}{\xi} \cdot \tilde{c}(\zeta)_g$ for all $g \in G$. Putting $s = s'\xi^{-1}$ we still have $s^k = t$ and $\frac{g(s)}{s} = \tilde{c}(\zeta)_g$ for all $g \in G$, i.e., $\frac{(\baru,\bare)(s)}{s} =  \zeta^e$, which is equivalent to $(\baru,\bare)(s) = \zeta^e \cdot s$. Hence, as mentioned above, $s$ generates $L$ over $K(\mu_k)$. \qed

\subsection{Applications to the Lubin-Tate tower} In this section we consider the fields $K_n$ and $\tK_n$ constructed in Section~\ref{defining_Kn}, but we assume throughout that $\cO = \Zp$ and $h=2$. The following is Result~\ref{result2i} (i) of the introduction.

\begin{cor}\label{prop:weak_Kummer}
The extension $K_n \mid K_0$ is a Kummer extension, i.e., there is $t_n \in K_0$ such that $K_n = K_0(\mu_{p^n},\sqrt[p^n]{t_n})$. The same is true for $\tK_n \mid \tK_0$.
\end{cor}

\Pf By Proposition~\ref{prop:gal-groups-knprime}, $K_n \mid K_0$ and $\tK_n \mid \tK_0$ are both Galois extensions with Galois group $(\Z/p^n\Z)^\times \ltimes \Z/p^n\Z$. Let $\bQp$ be the completion of the maximal unramified extension of $\Qp$. It has been shown in \cite[Cor.~3.4]{StrauchGeometricallyConnectedComponents} that the field $K'_n$ contains $\bQp(\mu_{p^n})$, which is a Lubin-Tate extension for the multiplicative formal group over $\bQp$. Since $K'_n \subseteq \tK_n \subseteq K_n$, both fields $K_n$ and $\tK_n$ contain $\mu_{p^n}$. By loc.\ cit., the action of 

\[\Gal(K'_n \midc K'_0) \simeq \left\{g_{a,b,d} := \left(\begin{array}{cc} a & b \\ 0 & d \end{array}\right) \; \Big| \; a,d \in (\Z/p^n\Z)^\times \,,\; b \in \Z/p^n\Z\right\} \;,\]

\vskip8pt

cf. \ref{prop:galois-group-k}, on $\mu_{p^n}$ is given by $g_{a,b,d}.\zeta  = \zeta^{ad}$, $\zeta \in \mu_{p^n}$. The subgroup $\Gal(\tK_n \mid \tK_0) = \Gal(K_n \mid K_0)$ of $\Gal(K'_n \midc K'_0)$ consists precisely of those $g_{a,b,d}$ with $d=1$. Therefore, elements $g_{a,b,1} \in \Gal(\tK_n \mid \tK_0) = \Gal(K_n \mid K_0)$ act on $\zeta \in \mu_{p^n}$ as $g_{a,b,1}.\zeta  = \zeta^a$. We are thus in the situation of Proposition~\ref{prop_Kummer} from which our assertion follows.
\qed

\vskip8pt

\begin{para} We now turn to the question whether the tower $K_\bullet$ is a Kummer tower. As mentioned in the introduction, for $K_\bullet$ to be a Kummer means that there is a $t \in K_0$ such that for every $n \ge 0$ one has $K_n = K_0(\mu_{p^n},\sqrt[p^n]{t})$. Our methods, however, are such that we can only investigate this question under the restriction that $t$ lies in the field $\tK_0$ of which $K_0$ is the $p$-adic completion. The point is that $\tK_n \midc K'_0$ are Galois extensions whereas $K_0 \mid K'_0$ is not, cf. \ref{GaloisNonGalois}, and our methods are tied to the fact that $\tK_n \midc K'_0$ is a Galois extension.\footnote{A more sophisticated approach might possibly give the stronger result about the non-existence of such a $t$ in $K_0$ (and not only in $\tK_0$), but we do not know how to do this.}

\vskip8pt

Recall the field $K'_{\ell,u} \subseteq K'_\ell$ which is the maximal unramified subextension of $K'_\ell \mid K'_0$, cf. \ref{defining_Kn}. We have, by definition, $\tK_0 = \bigcup_{m \ge 0} K'_{\ell,u}$, and $\tK_n = \tK_0 K'_n$. We recall from \ref{prop:gal-groups-knprime} (ii) that the universal Drinfeld basis induces an isomorphism

\[G_{n,\ell} := \Gal(\tK_n \midc K'_{\ell,u}) \simeq \left\{\begin{pmatrix} a&b\\0&d \end{pmatrix}\ \middle|\ a\in(\Z / p^n\Z)^\times \,,\,\ b\in\Z / p^n\Z \,,\,\ d\in 1+ p^\ell\Zp \right\} \;,\]

\vskip8pt

where $1+p^\ell\Zp$ is to be interpreted as $\Z_p^\times$ when $\ell=0$. In the remainder of this paper we use the notation $g_{a,b,d} = \begin{pmatrix} a&b\\0&d \end{pmatrix}$, as in the proof of \ref{prop:weak_Kummer}.
\end{para}

\vskip8pt

\begin{prop}\label{cohomologyGnell} Suppose $p>2$. The group $G_{n,\ell} = \Gal(\tK_n \midc K'_{\ell,u})$ acts on $\mu_{p^n} \subseteq K'_n \subseteq \tK_n$ by $g_{a,b,d} \cdot \zeta = \zeta^{ad}$ for $\zeta\in\mu_{p^n}$. Furthermore, 

\[H^1(G_{n,\ell},\mu_{p^n}) \simeq \{\zeta \in \mu_{p^n} \midc \forall d\in\bbZ_p^\times: \zeta^{d^2} = \zeta\ \} \;.\]

\vskip8pt

If $p>3$, then $H^1(G_{n,\ell},\mu_{p^n}) = \mu_{p^\ell}$.
\end{prop}

\Pf As we have already recalled in the proof of Proposition~\ref{prop:weak_Kummer}, the first assertion about the action of $G_{n,\ell}$ on $\mu_{p^n}$ is \cite[Cor.~3.4]{StrauchGeometricallyConnectedComponents}. Let

\[G_n := \Gal(\tK_n \midc \tK_0) = \{g_{a,b,d} \in G_{n,\ell} \ |\ d=1\} \simeq (\Z/p^n\Z)^\times \ltimes \Z/p^n\Z \;;\]

\vskip8pt

cf. \ref{prop:gal-groups-knprime} (iii). Then $G_n$ is a normal subgroup of $G_{n,\ell}$. Put $D = G_{n,\ell}/G_n \simeq 1+ p^\ell\Zp$. Since $\mu_{p^n}^{G_n} = \{1\}$, it follows from the inflation-restriction sequence \cite[6.8.3]{Weibel} that

\[H^1(G_{n,\ell},\mu_{p^n}) \simeq H^0(D, H^1(G_n,\mu_{p^n})) \;.\]

\vskip8pt

By Proposition~\ref{prop:H1_semidir}, we have $H^1(G_n,\mu_{p^n}) \simeq \mu_{p^n}$, where $\zeta\in\mu_{p^n}$ corresponds to the class of $c(\zeta)_{g_{a,b,1}} = \zeta^b$. The projection $G_{n,\ell} \twoheadrightarrow D$ has the section $D \ra G_{n,\ell}$, $d \mapsto \tilde{d} := g_{1,0,d}$. Then $\tilde{d}^{-1} g_{a,b,1} \tilde{d} = g_{a,bd,1}$, and so

\[(\tilde{d} \cdot c(\zeta))_{g_{a,b,1}} = \tilde{d} \cdot c(\zeta)_{\tilde{d}^{-1} g_{a,b,1} \tilde{d}} = \tilde{d} \cdot c(\zeta)_{g_{a,bd,1}} = \tilde{d} \cdot \zeta^{bd} = \zeta^{bd^2} \;.\]

\vskip8pt

Thus the cocycle is fixed by $D$ if and only if $\zeta^{d^2} = \zeta$ for all $d \in 1+p^\ell\Zp$. If $p>3$, then there exists a $d \in 1+p^\ell\Zp$ such that $p^{\ell+1}\nmid d^2-1$. Then $\zeta^{d^2} = \zeta$ if and only if $\zeta \in \mu_{p^\ell}$.
\qed

\begin{lemma}\label{lemma-invariants} The action of $G_{n,\ell} = \Gal(\tK_n \mid K'_{\ell,u})$ on $\tK_n$ extends by continuity to an action on $K_n$, and $H^0(G_{n,\ell}, K_n) = K'_{\ell,u}$. 
\end{lemma}

\Pf As the Galois automorphisms in $G_{n,\ell}$ are $p$-adically continuous, they extend uniquely to the $p$-adic completion $K_n$ of $\tK_n$. Consider the tautological exact sequence

\[1 \lra G_n = \Gal(\tK_n \mid \tK_0) \lra G_{n,\ell} = \Gal(\tK_n \mid K'_{\ell,u}) \lra D := \Gal(\tK_0 \mid K'_{\ell,u}) \lra 1 \;.\]

\vskip8pt

Recall that the canonical map $\Gal(K_n \mid K_0) \ra \Gal(\tK_n \mid \tK_0)$ is an isomorphism, by Proposition~\ref{prop:gal-groups-knprime} (iv). We therefore have 

\[H^0(G_{n,\ell},K_n) = H^0\Big(D,H^0(G_n,K_n)\Big) = H^0\Big(D,K_0\Big) \;.\]

\vskip8pt

Recall that $K_0$ is, by definition, the $p$-adic completion of $\tK_0$, cf \ref{defining_Kn}. Recall also that we defined $\tk_0 = \bigcup_n k'_{n,\sep}$ to be the residue field of $\tK_0$ (which is also the residue field of $K_0$), cf. \ref{defining_Kn}. The field $\tk_0$ is a Galois extension of $k'_{\ell,\sep}$ whose Galois group is canonically isomorphic to $D \simeq 1+p^\ell\Zp$. If we put $K''_0 = H^0(D,K_0)$, then group $D$ acts trivially on $K''_0$, and it acts therefore trivially as on its residue field, which must then be $k'_{\ell,\sep}$. The field $K''_0$ is therefore a discretely valued complete subfield of $K_0$ with residue field $k'_0$ and must then be equal to $K'_0$.
\qed

\begin{comment}

\begin{cor}\label{corollary_pn_powers}
If $p>3$, then 

\[(K_n^\times)^{p^n} \cap (K'_0)^\times = ((K'_0)^\times)^{p^n} \;.\]

\vskip8pt

\end{cor}

\Pf We apply $H^0(G_{n,\ell},-)$ to the short exact sequence $1 \ra \mu_{p^n} \ra K_n^\times \ra (K_n^\times)^{p^n} \ra 1$. Using \ref{lemma-invariants}, we obtain the exact sequence

$$(K'_0)^\times \xrightarrow{\;\;\;(\cdot)^{p^n}\;\;} (K_n^\times)^{p^n} \cap (K'_0)^\times \lra H^1(G_{n,\ell},\mu_{p^n}) \;.$$

\vskip8pt

By \ref{cohomologytildeGn} the cohomology group on the right vanishes, which proves our assertion. 
\qed

\begin{cor} There is no element $t \in K'_0$ such that $K_n = K_0(\mu_{p^n}, \sqrt[p^n]{t})$ for all $n \geq 0$.
\end{cor}

\Pf
Suppose there exists such a $t$. As $t$ is a $p^n$-th power in $K_n$, it belongs to $(K_n^\times)^{p^n} \cap (K'_0)^\times$, which is equal to $((K'_0)^\times)^{p^n}$ by \ref{corollary_pn_powers}. So the element $t^{1/p^n}$ belongs to $K_0(\mu_{p^n})$, which is properly contained in $K_n$ (which is easily seen by comparing degrees). Thus $t^{1/p^n}$ cannot generate $K_n$ over $K_0(\mu_{p^n})$.
\qed

\end{comment}

\begin{prop} Suppose $p>3$. There is no $t \in \tK_0$ such that for all sufiiciently large  $n \gg 0$ one has $K_n = K_0(\mu_{p^n}, \sqrt[p^n]{t})$.
\end{prop}

\Pf Suppose on the contrary that such a $t \in \tK_0$ exists. Then it is contained in some subfield $K_{\ell, u} \subseteq \tK_0$. Choose $n > \ell$. By \ref{lemma-invariants}, the action of $G_{n,\ell} = \Gal(\tK_n \mid K'_{\ell,u})$ on $\tK_n$ extends by continuity to an action on $K_n$. We can thus consider the 1-cocycle

\[G_{n,\ell} \ra \mu_{p^n} \;,\;\; s \mapsto s(t^{1/p^n})/t^{1/p^n} \;.\]

\vskip8pt

Because $H^1(G_{n,\ell},\mu_{p^n}) = \mu_{p^\ell}$, cf. \ref{cohomologyGnell},  the map $s \mapsto (s(t^{1/p^n})/t^{1/p^n})^{p^\ell}$ is a coboundary, so there exists a $\zeta \in \mu_{p^n}$ such that 

\[s(t^{1/p^{n-\ell}})/t^{1/p^{n-\ell}} = s(\zeta)/\zeta \;.\]

\vskip8pt

If $s \in \Gal(\tK_n \mid K_{\ell, u}(\mu_{p^n}))$, then we further have $s(t^{1/p^{n-\ell}})/t^{1/p^{n-\ell}} = 1$. Using \ref{lemma-invariants} again, we  conclude that $t^{1/p^{n-\ell}} \in K_{\ell, u}(\mu_{p^n})$. It then follows that $[K'_{\ell,u}(\mu_{p^n}, t^{1/p^n}) : K'_{\ell, u}(\mu_{p^n})] \leq p^\ell$, and thus

\[[K_0(\mu_{p^n}, t^{1/p^n}) : K_0] \leq [K_{\ell,u}(\mu_{p^n}, t^{1/p^n}) : K_{\ell,u}] \leq (p - 1) p^{n-1} p^\ell \;. \]

\vskip8pt

But $[K_n : K_0] = (p - 1) p^{n-1} p^n$. Thus there cannot exist such an element $t \in \tK_0$.
\qed

\bibliographystyle{plain}
\bibliography{RefList}

\def\cprime{$'$}
\begin{thebibliography}{10}

\bibitem{Bourbaki}
Nicolas Bourbaki.
\newblock {\em Elements of mathematics. {C}ommutative algebra}.
\newblock Hermann, Paris; Addison-Wesley Publishing Co., Reading, Mass., 1972.
\newblock Translated from the French.

\bibitem{Drinfeld}
V.~G. Drinfel{\cprime}d.
\newblock Elliptic modules.
\newblock {\em Math. USSR-Sb.}, 23(4):561--592, 1974.

\bibitem{FriedJarden}
Michael~D. Fried and Moshe Jarden.
\newblock {\em Field arithmetic}, volume~11 of {\em Ergebnisse der Mathematik
  und ihrer Grenzgebiete. 3. Folge. A Series of Modern Surveys in Mathematics
  [Results in Mathematics and Related Areas. 3rd Series. A Series of Modern
  Surveys in Mathematics]}.
\newblock Springer-Verlag, Berlin, third edition, 2008.
\newblock Revised by Jarden.

\bibitem{Kohlhaase}
Jan Kohlhaase.
\newblock Admissible {$\phi$}-modules and {$p$}-adic unitary representations.
\newblock {\em Math. Z.}, 270(3-4):839--869, 2012.

\bibitem{LubinTate}
Jonathan Lubin and John Tate.
\newblock Formal moduli for one-parameter formal {L}ie groups.
\newblock {\em Bull. Soc. Math. France}, 94:49--59, 1966.

\bibitem{Scholl}
Anthony~J. Scholl.
\newblock Higher fields of norms and {$(\phi,\Gamma)$}-modules.
\newblock {\em Doc. Math.}, Extra Vol.:685--709 (electronic), 2006.

\bibitem{Serre_LCFT}
J.-P. Serre.
\newblock Local class field theory.
\newblock In {\em Algebraic {N}umber {T}heory ({P}roc. {I}nstructional {C}onf.,
  {B}righton, 1965)}, pages 128--161. Thompson, Washington, D.C., 1967.

\bibitem{Serre_LF}
Jean-Pierre Serre.
\newblock {\em Local fields}, volume~67 of {\em Graduate Texts in Mathematics}.
\newblock Springer-Verlag, New York-Berlin, 1979.
\newblock Translated from the French by Marvin Jay Greenberg.

\bibitem{StrauchDeformationSpaces}
Matthias Strauch.
\newblock Deformation spaces of one-dimensional formal modules and their
  cohomology.
\newblock {\em Adv. Math.}, 217(3):889--951, 2008.

\bibitem{StrauchGeometricallyConnectedComponents}
Matthias Strauch.
\newblock Geometrically connected components of {L}ubin-{T}ate deformation
  spaces with level structures.
\newblock {\em Pure Appl. Math. Q.}, 4(4, Special Issue: In honor of
  Jean-Pierre Serre. Part 1):1215--1232, 2008.

\bibitem{StrauchGaloisOnTorsion}
Matthias Strauch.
\newblock Galois actions on torsion points of one-dimensional formal modules.
\newblock {\em J. Number Theory}, 130(3):528--533, 2010.

\bibitem{Weibel}
Charles~A. Weibel.
\newblock {\em An introduction to homological algebra}, volume~38 of {\em
  Cambridge Studies in Advanced Mathematics}.
\newblock Cambridge University Press, Cambridge, 1994.

\end{thebibliography}

\end{document}